\documentclass[a4paper,twoside]{article}

\usepackage{amsmath,amsthm,amsfonts,latexsym,amscd,amssymb,enumerate}

\usepackage[hypertex,backref]{hyperref}

\swapnumbers
\theoremstyle{plain}
\newtheorem{theorem}{Theorem}[section]
\newtheorem{lemma}[theorem]{Lemma}
\newtheorem{corollary}[theorem]{Corollary}
\newtheorem{proposition}[theorem]{Proposition}

\theoremstyle{definition}
\newtheorem{definition}[theorem]{Definition}
\newtheorem{example}[theorem]{Example}
\newtheorem{notation}[theorem]{Notation}
\newtheorem{convention}[theorem]{Convention}

\newtheorem{set-up}[theorem]{Geometric set-up}
\newtheorem{remark}[theorem]{Remark}

\newcommand{\reals}{\mathbb{R}}
\newcommand{\complexs}{\mathbb{C}}
\newcommand{\naturals}{\mathbb{N}}
\newcommand{\integers}{\mathbb{Z}}
\newcommand{\rationals}{\mathbb{Q}}

\newcommand{\boundary}[1]{\partial#1}

\newcommand{\abs}[1]{\left\lvert#1\right\rvert} 

\newcommand{\tensor}{\otimes}
\newcommand{\into}{\hookrightarrow}
\newcommand{\onto}{\twoheadrightarrow}
\newcommand{\iso}{\cong}
\newcommand{\congruent}{\equiv}
\newcommand{\disjointunion}{\amalg}

\newcommand{\semiProd}{\rtimes}
\newcommand{\semiprod}{\semiProd}

\DeclareMathOperator{\Diffeo}{Diffeo}
\DeclareMathOperator{\im}{im}      

\DeclareMathOperator{\spin}{spin}
\DeclareMathOperator{\charPos}{charPos^{\spin}}
\DeclareMathOperator{\Pos}{Pos} 

 \DeclareMathOperator{\tr}{tr}

\DeclareMathOperator{\ind}{ind}
\DeclareMathOperator{\Ind}{Ind}

\DeclareMathOperator{\scal}{scal}  

\newcommand{\forget}[1]{}

\def  \nuint {\raise10pt\hbox{$\nu$}\kern-6pt\int}

\newcommand\Tr{\operatorname{Tr}}

\newcommand\R{\mathcal R}
\newcommand\V{\mathcal V}
\newcommand\C{\mathcal C}

\def \Sp {{\cal S}}
\newcommand\B{\mathcal B}

\newcommand\D{\mathcal D}
\newcommand\Di{D\kern-6pt/}
\newcommand\cDi{{\mathcal D}\kern-6pt/}
\newcommand\spi{S\kern-6pt/}
\newcommand \cspi{\Sp\kern-6pt/}

\newcommand\CC{\mathbb C}

\def \cal {\mathcal}
\def \C {{\cal C}}

\newcommand\ZZ{\mathbb Z}

\newcommand\pa{\partial}

{\catcode`@=11\global\let\c@equation=\c@theorem}

\allowdisplaybreaks[2]

\begin{document}
\pagestyle{myheadings}
\markboth{Paolo Piazza and Thomas Schick}{Groups with torsion, bordism and rho-invariants}


\title{Groups with torsion, bordism and rho-invariants}

\author{Paolo Piazza and Thomas Schick}
\maketitle

\begin{abstract}
Let $\Gamma$ be a discrete group, and let $M$ be a
closed spin manifold of dimension $m>3$ with $\pi_1(M)=\Gamma$. We assume that $M$ admits a 
Riemannian metric of positive scalar curvature. 
We discuss how to use the $L^2$-rho invariant $\rho_{(2)}$ and 
the delocalized eta invariant $\eta_{<g>}$  associated to the Dirac operator on $M$ in order to
get information about the space of metrics with positive scalar
curvature.

In particular we prove that, if  $\Gamma$ contains
  {\it torsion} and  $m\equiv 3 \pmod 4$ then
  $M$ admits infinitely many different bordism classes
  of metrics with positive scalar curvature. This implies that there exist
  infinitely many concordance classes; we show that this is true 
  even up to
  diffeomorphism.

  If $\Gamma$ has certain special properties, e.g.~if it contains polynomially
  growing conjugacy classes of
  finite order elements, then  we obtain more refined information about the ``size''
  of the space of metric of positive scalar curvature, and these
  results also apply if the dimension is congruent to $1$ mod $4$. For
  example, if $\dim(M)\equiv 1 \pmod 4$ and $\Gamma$ contains a central element of odd order, then the
  moduli space of metrics of positive scalar curvature (modulo the action of
  the diffeomorphism group) has infinitely many components, if it is not empty.

  Some of our invariants are the delocalized eta-invariants introduced by
  John Lott. These invariants are defined by certain integrals whose
  convergence is not clear in general, and we show, in effect, that 
  examples exist where this integral definitely does not converge, thus answering a
  question of Lott.
 
  We also discuss the possible values of the rho-invariants of the
  Dirac operator and show that there are certain global restrictions
  (provided that the scalar curvature is positive).
 
\end{abstract}

\tableofcontents

\section{Introduction and main results}\label{sec:intro}

Let  $M$ be a closed smooth manifold with fundamental group $\Gamma$ and universal cover
 $\overline{M}$. In this paper, we are concerned mainly with the
set  $\R^+ (M)$ of metrics of positive scalar curvature  on $M$
(this is in fact a topological space).

There is of course a preliminary question, namely whether this space is non-empty.
It  is known that there are powerful obstructions to the existence of positive scalar curvature
($\equiv$ PSC) metrics, the most successful 
  being  the one implied by the Lichnerowicz formula: on a spin manifold with positive scalar
  curvature, the spin Dirac operator twisted by the Mishchenko line bundle
  $\V:=\overline{M}\times_\Gamma C^*_r\Gamma$ is invertible.
  In this paper we shall leave the existence problem aside and assume
  that there exist a metric with positive scalar curvature. We shall
  instead concentrate on the classification question;
if one such metric exists, how many can we put on $M$ that are {\it distinct}?
We need to  clarify what we mean by distinct.
There are   three ways for distinguishing  two metrics of positive scalar curvature
 $g_1$ and $g_2$ on $M$

The first one is to say that $g_1$ and $g_2$  are not
path-connected in $\R^+(M)$. Thus, in this case, we are interested in $\pi_0({\cal
R}^+(M))$, the set of arcwise connected components of ${\cal
R}^+(M)$.

The second  way for distinguishing two PSC metrics employs the
notion of {\it concordance}: $g_1$ and $g_2$ are concordant if
there exists a metric of PSC on $M\times [0,1]$ extending $g_1$ on
$M\times \{0\}$,
 $g_2$ on $M\times \{1\}$ and of product-type near the
boundary. The set of concordance classes of PSC
metrics on  $M$ is denoted by $\widetilde{\pi}_0( \R^+(M))$.

 \begin{convention}
    Throughout the paper, whenever we work with a Riemannian metric on a
    manifold with boundary, we assume that the metric has product structure
    near the boundary. Observe, in particular, that the restriction to the
    boundary of such a metric has positive scalar curvature, if the original
    one has positive scalar curvature.
  \end{convention}

The third and more subtle way for distinguishing two PSC metrics
$g_1$, $g_2$ on a spin manifold $M$  employs the notion of {\it
bordism}. 

  \begin{definition}\label{def:of_metric bordism}
    Let $M$ be a closed spin manifold with fundamental group  $\Gamma$. Two
    metrics $g_1$ and $g_2$ of positive scalar curvature on $M$ are
    \emph{$\pi_1$-spin bordant} if there is a compact spin manifold $W$ with
    positive scalar curvature metric $g$ and with boundary $\boundary W=
    (M,g_1)\disjointunion (-M,g_2)$, which admits a $\Gamma$-covering
       $\overline W$ whose boundary is the union of the universal
    coverings of the two boundary
    components. 
     
    Note that this notion has an evident extension to metrics on possibly
    different spin manifolds.
     \end{definition}

   It is obvious that if two metrics $(g_1,g_2)$ are  \emph{concordant} then they are 
   in particular bordant since we can choose $W=[0,1]\times M$ as the underlying manifold of
    the bordism. (On the other hand there are examples of
non-concordant metrics that are bordant, see \cite[page 329]{lawson89:_spin}.)
    It is also rather clear that two metrics which lie in the
    same path component of the space of all metrics of positive scalar
    curvature on a given manifold $M$ are concordant
    and, therefore, bordant.
    Summarizing, as far as the problem of distinguishing metrics of
positive scalar curvature is concerned, we have: $$\text{ non-bordant }
\,\,\Rightarrow\,\, \text{ non-concordant}\,\, \Rightarrow\,\, \text{ non-pathconnected. } 
    $$
     
     In this paper we shall use the $L^2$-rho invariant $\rho_{(2)}$ 
    of Cheeger-Gromov
     and the delocalized
     eta invariant $\eta_{<g>}$ of Lott for
     the spin Dirac operator associated to $(M,g)$ in order to distinguish 
     non-bordant metrics of
     positive scalar curvature.  Fundamental to our analysis will be the bordism invariance of 
           $\rho_{(2)}$ and  $\eta_{<g>}$, the long exact sequence of bordism  groups due to Stephan Stolz
           and some fundamental examples due to Botvinnik and Gilkey.

%
           In order to apply our methods, pioneered by Botvinnik and Gilkey in \cite{BoGi},
  we shall need to assume that $\Gamma \equiv \pi_1 (M)$ is not torsion-free: indeed
     if $\Gamma$ is torsion free and satisfies the Baum-Connes conjecture for the maximal
     group $C^*$-algebra then, because of the PSC assumption,
      these invariants are identically zero, as we have 
     proved in \cite{piazza:_bordis} \footnote{for the delocalized eta invariant it suffices
     to assume that $\Gamma$ satisfies the Baum-Connes conjecture for the {\it reduced} $C^*$-algebra};
     moreover, there are no known examples of torsion-free groups for which these invariants
     are non-zero.

As an example of the results we shall establish, we anticipate one of our main theorems:

\begin{theorem}\label{theo:L2-botvinnik-gilkey}
  Assume that $M$ is a spin-manifold of dimension $4k+3$, where $k>0$.
Assume that $g$ is a metric with positive scalar curvature
  on $M$, and that the fundamental group $\Gamma$ of $M$  contains torsion.
Then $M$ admits infinitely
  many different $\Gamma$-bordism classes of metric with $\scal>0$; they are
  distinguished by $\rho_{(2)}$. These infinitely many bordism classes
  remain distinct even after we mod out the action of the 
  group of diffeomorphisms of $M$\footnote{for the precise meaning of
    this compare Theorem \ref{theo:infinitely-many}}.
  \end{theorem}

This theorem generalizes results of Botvinnik-Gilkey \cite{BoGi}, \cite{BoGi2};
 other generalizations of their results have appeared in Leichtnam-Piazza \cite{LPPSC}.

Under additional assumptions on the group $\Gamma$
 we shall be able to estimate the size of the set of equivalence classes of non-bordant metrics
           by proving that a free group of a certain rank acts freely on this  set. 
We want to single out one consequence of these results, which also apply in
dimensions $4k+1$:
\begin{corollary}\label{cor:4k+1}
  If $\Gamma$ is a discrete group which contains a central element of odd
  order, and if $M$ is a spin manifold with fundamental group $\Gamma$ and of
  dimension $4k+1$ which admits a metric with positive scalar curvature, then
  the moduli space of such metrics (modulo the action of the diffeomorphism
  group via pullback) has infinitely many components 
  \begin{equation}\label{cardinality}
  |\pi_0 (\mathcal{R}^+(M)/{\rm Diffeo}(M))|=\infty
  \end{equation}
\end{corollary}

To our knowledge, this is the first general result of this kind which applies in
dimension congruent to $1$ modulo $4$. 

Extensions to even dimensional manifolds with special fundamental groups
should be  possible
by combining the methods of the current paper with those of \cite{LPPSC}.

Note that, if the dimension is
congruent  to $3$ mod $4$, then it is  \emph{always} true
that \eqref{cardinality} holds; compare \cite[Theorem 7.7 of Chapter IV]{lawson89:_spin}.

\medskip

There is a very parallel story for the signature operator, where the condition
on positive scalar curvature is replaced by ``homotopy invariance'' ---stated
differently, 
one gets vanishing or classification results for the disjoint
union of one manifold with a homotopy equivalent second manifold.

For instance, if  $\Gamma$ is torsion-free and satisfies the Baum-Connes
           conjecture for the maximal $C^*$-algebra, we prove in \cite{piazza:_bordis}          that 
      $\rho_{(2)}$ and  $\eta_{<g>}$ vanish
on a manifold which is the disjoint union of two homotopy equivalent manifolds.
For $\eta_{(2)}$ this result is  originally due to Keswani \cite{Kes2}.

Similarly, the non-triviality result we give in Theorem   \ref{theo:L2-botvinnik-gilkey}
has a relative for the signature
operator; a result of Chang-Weinberger \cite{math.GT/0306247} which was actually
the motivation for our result, and also for its proof.

\begin{theorem}\label{theo:chang-weinberger}
(Chang-Weinberger) 
If $M$ is a compact oriented manifold of dimension $4k+3$, where $k>0$, such that
$\pi_1 (M)$ is not torsion-free, then there are infinitely many
manifolds that are homotopy equivalent to $M$ but not homeomorphic
to it.
\end{theorem}

Notice that, in particular, the structure set $\mathcal{S}(M)$ has infinite cardinality.
Chang and Weinberger ask in their paper about more precise results concerning
the ``size'' of the structure set if the fundamental group contains a lot of
torsion. 
In this  paper we  investigate the corresponding question  for the space
of metrics of positive scalar curvature, and use in particular the delocalized
eta-invariants of John Lott to get some positive results ---for precise
statements consult Theorem \ref{theo:infinitely-many-sharp} and Theorem
\ref{theo:diff-conj-class-1}. 
It should be possible, although technically more
difficult given that the boundary operator is not invertible, to extend the results
stated in Theorem \ref{theo:infinitely-many-sharp} and Theorem
\ref{theo:diff-conj-class-1} to the signature operator and the structure set of a fixed
manifold. We plan  to investigate this and further directions of research for
the signature operator  in
future work.

\medskip
Our results rely on the delocalized eta-invariats of Lott
\cite{MR2000k:58039}, applied in those situations where they are well defined
and one does not have any convergence problems (e.g.~for central group
elements). However, we give in Section \ref{sec:an-example-non} examples which
show that in general the convergence one hopes for does definitely not
occur, showing the limitations of this method.

Close relatives of the delocalized rho-invariants we consider are the
rho-invariants associated to virtual representations of dimension zero (we
explain the translation between the two points of view via ``Fourier
transform'' in Section \ref{sec:many-bordism-classes}). Given such a finite
dimensional virtual unitary
representation $[\lambda_1-\lambda_2]$, let $F_1$ and $F_2$ be the associated flat
vector bundles. Then the corresponding rho-invariant is simply
\begin{equation*}
   \rho_{\lambda_1-\lambda_2}(\Di)= \eta(\Di_{F_1})-\eta(\Di_{F_2}).
\end{equation*}
One might wonder what the possible values of the rho-invariants
are, if the group is not torsion-free.
For the signature operator and these APS-rho invariants a result
of this type has been proved by Guentner-Higson-Weinberger:

\begin{theorem}\label{theo:Guentner-Higson-Weinberger}
\cite[Theorem 7.1]{Hi-Gu-We}
Let  $M$ and $N$ be smooth,
closed, oriented, odd-dimensional manifolds which are oriented  homotopy equivalent.
Let $\Gamma:=\pi_1 (M)$
 and let $\lambda_1, \lambda_2:\Gamma \to U(d)$ be two unitary representations.
Let $R'$ be the smallest subring of $\rationals$ generated by $\ZZ$, $1/2$, and
  $o(\lambda_1(g))^{-1}$ and $o(\lambda_2(g))^{-1}$ for each
  $g\in\Gamma$. Here, $o(x)$ is the order of the group element $x\in
  U(d)$, and we set $(+\infty)^{-1}:=0$.
Then
$$ \rho_{\lambda_1-\lambda_2}(D^{{\rm sign}}_M)
- \rho_{\lambda_1-\lambda_2}(D^{{\rm sign}}_{N})\in R'.$$
\end{theorem}

We end the paper by proving 
the corresponding result in the positive scalar curvature context. It is no
surprise that we don't need to invert $2$, as is notoriously necessary in
L-theory contexts.

\begin{theorem}\label{prop:possible_values_of_APS_rho_in_general}
Let $M$ be a spin manifold with a Riemannian metric of positive scalar
curvature and let $\Di$
be the associated Dirac operator.
Let $\lambda_j$  be as in Theorem \ref{theo:Guentner-Higson-Weinberger} and
let $R$ be the smallest subring of $\rationals$ generated by $\integers$ and
$o(\lambda_1(G)^{-1}$, $o(\lambda_2(g)^{-1}$ for each $g\in\Gamma$ of finite order. Then
  \begin{equation*}
    \rho_{\lambda_1-\lambda_2}(\Di) \in R.
  \end{equation*}
\end{theorem}

\medskip
\noindent
{\bf Acknowledgements.} Part of this work was carried out during visits of the authors to
G{\"o}ttingen and  Rome funded by Ministero Istruzione
Universit{\`a} Ricerca, Italy (Cofin {\it Spazi di Moduli e Teorie
di Lie}) and  Graduiertenkolleg ``Gruppen und Geometrie''
(G{\"o}ttingen).

\section{Distinguishing metrics with positive scalar
curvature}

\label{sec:exampl-many-diff}

\subsection{Torsion in $\pi_1(M)$ and $\dim M\equiv 3 \pmod 4$}
\label{sec:many-bordism-classes}




  We will frequently talk about spin manifolds; we think of them as being
  manifolds with a given spin structure (so they correspond to \emph{oriented}
  manifolds, not to orientable manifolds).

 Let $X$ be any space. Then there is an exact sequence
of bordism groups due to Stephan Stolz, see \cite{stolz:_concor_psc}, \cite[p.~630]{Sto1}.

\begin{equation}\label{eq:pos_sequemnce}
\to \Omega^{\spin}_{n+1}(X)\xrightarrow{t} R^{\spin}_{n+1}(X) \xrightarrow{\delta} \Pos^{\spin}_n(X)\to
  \Omega_n^{\spin}(X)\to R^{\spin}_n(\pi_1(X))\to
\end{equation}

Recall the definition of the terms in the sequence:
\begin{definition} 
  a) First, $\Omega_*^{\spin}(X)$ is the singular spin bordism group of
  $X$, the set of closed spin manifolds with a reference map to
  $X$, modulo spin bordism.\\
  b) $\Pos^{\spin}_*(X)$ is the bordism group of
  spin manifolds with a given metric with $\scal>0$, with a reference
  map to $X$. A bordism in $\Pos^{\spin}_*(X)$ is a bordism of spin
  manifolds  as above, together with a metric with positive scalar
  curvature  which restricts to the given metrics at the boundary
  (with a product structure near the boundary).\\
  c) $R^{\spin}_*(X)$ is
  the set of compact spin-manifolds with boundary, where the boundary
  is equipped with a metric with positive scalar curvature, together
  with a reference map to $X$, modulo bordism. A bordism consists
  first of a bordism of the boundary, with a metric with $\scal>0$ as
  in the bordism relation for $\Pos^{\spin}_*$. Glue this bordism of
  the boundary to the original manifold with boundary, to get a closed
  manifold. The second part of a bordism in $R^{\spin}_*(X)$ is a zero
  bordism of the resulting closed manifold (together with an extension
  of the reference map).\\
  d) The maps in the exact sequence \eqref{eq:pos_sequemnce} are quite
  obvious:\\
    $R^{\spin}_{n+1}\to \Pos^{\spin}_n$ is given by taking the
  boundary, $\Pos^{\spin}_n\to\Omega^{\spin}_n$ by forgetting the
  Riemannian metric, and $\Omega^{\spin}_n\to R^{\spin}_n$ by
  understanding a closed manifold as a manifold with empty boundary,
  this boundary therefore having a metric with $\scal>0$. \\
  e) The sequence
  is exact by definition. It is also evident that this sequence is
  natural with respect to maps $X\to Y$, and each entry is a covariant
  functor with respect to such maps.
\end{definition}

\begin{definition}
  Let $M$ be a closed spin manifold with fundamental group
  $\Gamma$. Let $u\colon M\to B\Gamma$ be a classifying map for a
  universal covering (i.e.~an isomorphism on $\pi_1$). We set
  $\charPos(M,u):=\{[M,g,u]\in \Pos^{\spin}_*(B\Gamma)\}$. 
 These are all the
  different bordism classes of metrics with positive scalar curvature
  on $M$ (where bordisms are considered which respect the given map $u$,
  i.e.~include the data of the fundamental
  group, and where also the spin structure on $M$ is fixed once and
  for all). Note that $\charPos(M,u)$ is a subset of the group
  $\Pos^{\spin}_n(B\Gamma)$, but we can't expect that it is a subgroup.

In this section, we will study the set $\charPos(M,u)$; we denote the class
$[M,g,u]\in \charPos(M,u)$ by $[g]$.
\end{definition}


\begin{proposition}\label{prop:action_of_ker_on_charPos}
  Let $M$ be a spin manifold with fundamental group $\Gamma$, of
  dimension $m\ge 5$. If $$[g]\in \charPos(M,u)\subset
  \Pos_m^{\spin}(B\Gamma)\,\text{ and }\,x\in\ker(\Pos^{\spin}_m(B\Gamma)\to
  \Omega_m^{\spin}(B\Gamma))$$then $x+[g]\in \charPos(M,u)$. 
  The action of $\ker(\Pos^{\spin}_m(B\Gamma)\to
  \Omega_m^{\spin}(B\Gamma))$ on $\charPos(M,u)$ is free and transitive.
\end{proposition}
\begin{proof}
  The statement is a consequence of the surgery result of
  Gromov-Lawson, Schoen-Yau, Gajer (compare \cite[Lemma 3.1]{BoGi}). 
  Since the underlying manifold
  $(X,f\colon X\to B\Gamma)$ of
  $x$ is zero bordant in $\Omega^{\spin}_m(B\Gamma)$, the sum of $(M,u)$
  and $(X,f)$ is bordant to $(M,u)$ in
  $\Omega_m^{\spin}(B\Gamma)$. By assumption, on this sum we have a
  metric with $\scal>0$. This metric can, by the surgery result, be
  extended over a suitable modification of the interior of this bordism (we have to
  make it sufficiently connected), to yield some metric with $\scal
  >0$ on the other end, i.e.~on $(M,u)$. If we perform two such
  constructions, we can glue the resulting bordisms (with their metric
  and reference map to $B\Gamma)$ along the boundary $(M+X,u+f)$ to see that
  the bordism class of the resulting metric is well defined.

  Since the action comes from addition in $\Pos^{\spin}_n(B\Gamma)$,
  the statement about 
  freeness 
  follows immediately. In order to prove transititivity we simply observe
   that any two objects $[g_1],[g_2]\in\charPos(M,u)$ are
  mapped to the same element of  $\Omega^{\spin}_n(B\Gamma)$, so that
  their difference belongs to $\ker(\Pos^{\spin}_m(B\Gamma)\to
  \Omega_m^{\spin}(B\Gamma))$. Thus 
  $$[g_2]=[g_1]+([g_2]-[g_1])=[g_1] + x\,,\,\,\text{with}\,\, x\in\ker(\Pos^{\spin}_m(B\Gamma)\to
  \Omega_m^{\spin}(B\Gamma))
  $$
  and we are done.
%
\end{proof}


Now, we want to introduce invariants on $\Pos^{\spin}_m(B\Gamma)$ and
$\charPos(M,u)$ which can be used to distinguish elements in these
sets.

\begin{definition}\label{def:rho_2}
  Let $(M,g)$ be a spin-manifold with Riemannian metric $g$ and with
  reference map $u\colon M\to B\Gamma$. Let $\overline M$ be the
  $\Gamma$-covering classified by $u$ (if $u$ is a
  $\pi_1$-isomorphism, then $\overline M$ is a universal covering of
  $M$). Define
  \begin{equation*}
    \rho_{(2)}(M,g,u):=  \eta_{(2)}(\overline D)- \eta(D) ,
  \end{equation*}
  where $D$ is the spin Dirac operator on $M$ and $\overline D$ its
  lift to $\overline M$. For details on the eta and the $L^2$-eta-invariant, compare, for example,  \cite{piazza:_bordis}.

  Fix an element $h\in\Gamma$ such that its conjugacy class $<h>$ has
  polynomial growth (inside $\Gamma$ with its word metric). If the scalar
  curvature of $(M,g)$ is strictly positive, then the Dirac operator of $M$
  (and the Dirac operator twisted with any flat bundle) is
  invertible. Consequently, the 
  delocalized eta invariant of Lott, denoted $\eta_{<h>}(\overline{D})$, is defined (compare
  \cite{MR2000k:58039} \cite{delocalized-erratum}; see also  \cite[Section 13.1]{piazza:_bordis}). More precisely,
  \begin{equation*}
   \eta_{<h>}(\overline{D}) = \frac{1}{\sqrt{\pi}}\int_0^\infty \sum_{\gamma\in<h>}
    \left(\int_{\mathcal{F}} \tr_x 
    k_t(x,\gamma x)\,dx\right)\,\frac{dt}{\sqrt{t}}, 
  \end{equation*}
  here $k_t(x,y)$ is the integral kernel of the operator $\overline D
  e^{-t\overline 
  D^2}$ on the covering $\overline{M}:=u^* E\Gamma$ and  $\mathcal{F}$ is a fundamental
  domain for this covering.  Note that it is a highly non-trivial
  fact that this sum and integral converge; it is proved for invertible
  $\overline D$ and groups  $\Gamma$ of polynomial growth in  \cite{MR2000k:58039}; 
  we observed in 
  \cite{piazza:_bordis} that one can take arbitrary groups, provided the 
  conjugacy class is of polynomial growth. Some information about conjugacy classes of polynomial growth 
can
be found in \cite{Wiethaup}. We give an example
  where the expression does not converge in Section \ref{sec:an-example-non}.
Notice that the same formula, if $h=1$, defines
  $ \eta_{(2)}(\overline D)$.
\end{definition}

\begin{notation}\label{notat:delocalized-rho} 
If $h\ne 1$ we shall  set
$ \rho_{<h>}(M,g,u):=\eta_{<h>}(\overline D)$.
\end{notation}

\begin{proposition}\label{prop:rho_2_as_homom}
The invariants  $\rho_{(2)}$ and $\rho_{<h>}$ of Definition \ref{def:rho_2} define  homomorphisms
  \begin{equation*}
    \rho_{(2)}\colon \Pos_*^{\spin}(B\Gamma)\to
    \reals,\qquad\rho_{<h>}\colon \Pos_*^{\spin}(B\Gamma)\to
    \complexs.
  \end{equation*}
\end{proposition}
\begin{proof}
  The group structure in $\Pos_m^{\spin}(B\Gamma)$ is given by disjoint
  union, and $\rho_{(2)}$ as well as $\rho_{<h>}$ are additive under disjoint
  union. 
  We only have to show that the invariant vanishes for a manifold
  representing $0$ in $\Pos_m^{\spin}(B\Gamma)$.
  Similar arguments have already been
  used in many places, e.g. \cite[Theorem 1.1]{BoGi}
and \cite[Proposition 4.1]{LPPSC}.
Let $[M,g,u]=0$ in $\Pos_m^{\spin}(B\Gamma)$; then there exists a spin Riemannian manifold
$(W,G)$ together with a classifying map $U:W\to B\Gamma$ such that $\pa W=M$,
$U|_{\pa W}=u$, $\scal(G)>0$, $G|_{\pa W}=g$.
Let $\D_W$ be the Mishchenko-Fomenko spin Dirac operator associated to 
$(W,G,U)$; let $C^*\Gamma$ be the maximal group $C^*$-algebra.
Since $\scal(g)>0$ there is a well defined index class
$\Ind(\D_W)\in K_0 (C^* \Gamma)$; since $\scal(G)>0$ this index class is zero (the operator
$\D_W$ is in fact invertible). Let $\Ind_{[0]}(\D_W):= \Tr^{\rm alg}(\Ind (\D_W))$
with $\Tr^{\rm alg}: K_0 (C^* \Gamma)\to C^*\Gamma/\overline{[C^*\Gamma,C^*\Gamma]}$
the natural algebraic  trace. Then, the APS index theorem proved in \cite{LPMEMOIRS}
(see  \cite[Theorem 3.3]{piazza:_bordis} for a direct and elementary proof 
of the special case
used here) gives
\begin{equation}\label{vanish-index}
0=\Ind_{[0]} (\D_W)= (\int_W \widehat{A}(W,G))\cdot 1-\frac{1}{2}\eta_{[0]}(\D_{M})\,\in\,\,C^*\Gamma/\overline{[C^*\Gamma,C^*\Gamma]}
\,;\end{equation}
a similar identity holds in the abelianization $C^*_r\Gamma/\overline{[C^*_r\Gamma,C^*_r \Gamma]}$
of the reduced group $C^*$-algebra as well as in the abelianization 
of the Connes-Moscovici
algebra $\B^\infty_\Gamma$:
\begin{equation}\label{vanish-index-infty}
0=\Ind_{[0]} (\D_W^\infty)= (\int_W \widehat{A}(W,G))\cdot 1-\frac{1}{2}\eta_{[0]}(\D_{M}^\infty)\,\,\in\,\,
\B^\infty_\Gamma/\overline{[\B^\infty_\Gamma,\B^\infty_\Gamma]}
\,.\end{equation}
Let $h\not=0$ and let  $\tau_{<h>}\colon \complexs\Gamma\to\complexs$ be the trace 
defined by $  \sum_{g\in \Gamma}\lambda_g g\mapsto
  \sum_{g\in <h>} \lambda_g$. Because the conjugacy class $<h>$ has
  polynomial growth, we observed in \cite[Proposition 13.5]{piazza:_bordis}
  that $\tau_{<h>}$  extends by continuity to a  trace on
  $\B^\infty_\Gamma$. By  \cite[Formula (4.16)]{MR2000k:58039}, $$\tau_{<h>}(\eta_{[0]}(\D^\infty_{M}))
  =\eta_{<h>} (\overline{D})\equiv \rho_{<h>}(M,g,u)\,,$$ and since  $\tau_{<h>} (1)=0$,
  we finally see that by applying $\tau_{<h>} $ to (\ref{vanish-index-infty}) we get 
  $\rho_{<h>}(M,g,u)=0$ which is what we wanted to prove.\\
  Let $\tau:C^*\Gamma\to \CC$ be the trace induced by the trivial representation;
  let $\tau_\Gamma: C^*\Gamma\to \CC$ be the  canonical trace, i.e. the trace induced by 
 $ \sum_{g\in \Gamma}\lambda_g g\mapsto \lambda_1$. Obviously $\tau(1)=\tau_\Gamma (1)$.
Recall now that we have also proved in \cite{piazza:_bordis} that
$$\tau_{\Gamma}(\eta_{[0]}(\D_{M}))=\eta_{(2)}(\overline{D})\,,\quad \tau(\eta_{[0]}(\D_{M}))=\eta(D)\,;$$
we complete the proof of the Proposition by applying $\tau$ and $\tau_\Gamma$ to (\ref{vanish-index})
and 
subtracting.
  

\end{proof}

\begin{proposition}\label{prop:diffeo_and_rplus}
Let $\mathcal{R}^+(M)$ be the space of smooth metrics with positive
scalar curvature on $M$, with the usual $C^\infty$-topology. For the fixed
spin structure and classifying map $u\colon M\to B\Gamma=B\pi_1(M)$, we get an
obvious surjection $\mathcal{R}^+(M)\onto \charPos(M,u)$.
The composition 
  \begin{equation*}
    \mathcal{R}^+(M) \to \charPos(M,u) \into \Pos^{spin}_m(B\Gamma)
    \xrightarrow{\rho_{(2)}} \reals
  \end{equation*}
is constant on orbits of the action of the spin-structure preserving
  diffeomoprhism group $\Diffeo_{\#}(M)$ (which acts by pulling back
  the Riemannian metric).
 Moreover, it is locally constant, and therefore factors through the set of
  components of the moduli space $\pi_0(\mathcal{R}^+(M)/\Diffeo_{\#}(M))$.
  \end{proposition}
  \begin{proof}
   Let  $PSpin (TM)$ be  a $2$-fold covering  of $PSO_g (TM)\to M$ which is non-trivial on the
    fibers and wich  determines the chosen spin structure on $M$. 
    Equivalence of spin structures
    is understood as equivalence of such 2-fold coverings.
     Using  the fact that the inclusion $PSO_h (TM)\hookrightarrow PGL_+ (TM)$ 
     is a homotopy equivalence for each metric $h$ on $TM$,
     we can equivalently define a spin structure
     as a 2-fold covering of $PGL_+ (TM)$ which is non-trivial along the fibers
   of   $PGL_+ (TM)\to M$; this means, in particular, that the choice of a spin structure
   for one metric $g$ canonically determines a spin structure for any other
   metric $h$ ---compare  \cite[Chapter II, Sections 1 and 2]{lawson89:_spin}.
   Let $\Psi:M\to M$ be  a diffeomorphism; let $d\Psi: PGL_+ (TM)\to PGL_+
   (TM)$  be
   the induced diffeomorphism. Then $\Psi$ is spin structure preserving if the pullback
  $d\Psi^*(PSpin(TM))$  is equivalent to $PSpin(TM)$. Call the  corresponding isomorphism $\beta_{GL_+}$.
   Now, if we define the spinor bundle, $L^2$-spinors and the Dirac
  operator entirely in terms of the pullback structures, $\Psi$ induces
  a unitary equivalence, and consequently the eta invariant of $D$ and
  of the operator defined using the pulled back structure coincide. On
  the other hand, the isomorphism $\beta_{GL_+}$ induces an isomorphism 
  $\beta$ between the original spin structure and the pulled-back structure
  both seen as 2-fold coverings of  $PSO_{\Psi^*g} (TM)$; using $\beta$ we get 
   a unitary equivalence
  between the operator obtained via the pulled back structures and
  the Dirac operator for $\Psi^*g$ and the chosen fixed spin structure, so
  that their eta invariants coincide, as well. Taken together, $\eta(D_g)=\eta(D_{\Psi^*g})$.
      More or less the same applies to the construction of the $L^2$-eta invariant on the
    universal covering. In order to simplify the notation, let us
    denote by $P$ the chosen spin structure. We start with a given covering
    $\overline M\xrightarrow{\pi} M$ with given action of $\Gamma$ by deck transformations
    (obtained by pulling back $E\Gamma$ from $B\Gamma$ via the map
    $u\colon M\to B\Gamma$).  The spin structure and the metric on $\overline M$, 
    denoted $\overline g$ and $\overline P$, are 
    the ones pulled back
    from $M$ via $\pi$.
    We can then pull back everything, including the covering
    $\overline M$ via $\Psi$, and will obtain a $\Gamma$-covering
    $p\colon\Psi^*\overline M\to M$ with pullback $\Gamma$-action, spin structure, pullback
    metric etc. Then $\Psi$ will induce a unitary $\Gamma$-equivariant
    equivalence between $\overline D$ and the Dirac operator
    constructed entirely in terms of the pulled back structures, so
    the $L^2$-eta invariants of these two operators coincide.
On the other hand, we have the covering $\overline M$ itself, and
    the fixed spin structure. Since the universal covering is unique,
    we get a covering isomorphism $\gamma\colon \overline
    M\to\Psi^*\overline M$, covering the identity. It becomes an isometry if
   we use on $\overline M$ the lift of the metric $\Psi^*g$. 
   On $\Psi^*\overline M$ we 
   have used  the spin structure given by the pullback principle bundle
   $\overline\Psi^*\overline P$ with $\overline\Psi$ the obvious map
   $\Psi^* \overline M \to \overline M$  covering $\Psi$. 
  Since $ p\circ \gamma=\pi$, we get a map of
   principal bundles $\pi^*\Psi^* P\to \overline\Psi^*\overline P$. We now use
   the principal bundle isomorphism $P\to \Psi^*P$ which comes from the fact
   that $\Psi$ is spin structure preserving, to finally identify the spin principle
   bundle of $\overline M$ to the one of $\Psi^*\overline M$ via a map $\overline \gamma$
    covering 
   $\gamma$ and   the map $P\to \Psi^*P$ of principal bundles on $M$.
 Proceeding as for $M$ itself, we obtain
    a unitary equivalence between $\overline D_{\Psi^*g}$ and the
    operator obtained using the pullback structures.
    Summarizing: $\eta_{(2)}(\overline D _g)=  \eta_{(2)}(\overline D _{\Psi^*g})$.
   \end{proof}

\begin{remark}\label{remark:diffeo-and-delocalized} 
    It should be noted that the map $\gamma$ given above is not, in general, 
    $\Gamma$-equivariant, but we can choose $\gamma$ in such a way
    that for $x\in \overline M$ and $g\in \Gamma$, $\gamma(gx)=
    \alpha_\Psi(g)\gamma(x)$, where $\alpha_\Psi\colon\Gamma\to \Gamma$ is equal
    to the isomorphism $u_* \pi_1(\Psi) u_*^{-1}$. This is true
    because (by the universal property of $B\Gamma$ and $E\Gamma$),
    $\Psi^*u^*E\Gamma$ is isomorphic as $\Gamma$-principal bundle to
    $u^* (B\alpha_\Psi)^* E\Gamma$, since $u\circ \Psi$ and $B\alpha_\Psi\circ u$
    induce the same map on the fundamental group. Moreover, by
    \cite[Appendix B, p.~378]{lawson89:_spin},
    $(B\alpha_\Psi)^*E\Gamma$ is isomorphic as $\Gamma$-principal bundle to
    the associated bundle $\Gamma\times_{\alpha_\Psi} E\Gamma$, 
    and the required covering isomorphism
    \begin{equation*}
      E\Gamma\to
    \Gamma\times_{\alpha_\Psi} E\Gamma; x\mapsto [1,x]
  \end{equation*}
  has exactly the required equivariance property:
  $[1,gx]=[\alpha_\Psi(g),x]=\alpha_\Psi(g)[1,x]$, which is preserved when
  pulling back the whole covering isomorphism with $u$.
    
  Now,  as explained above,
     the map $\gamma$ induces maps which preserve all the structure
  which is present in the construction of the Dirac operators on
  $\overline M$ (using the lift of the metric $\Psi^*g$) and $\Psi^*\overline M$ (except for
  the group action). In particular, for the fiberwise trace we have
  $$\tr k_t(x,hx)
  = \tr \kappa_t(\gamma(x),\alpha_\Psi(h)\gamma(x))\,,$$
   where here $k_t(x,y)$ is
  the integral kernel of $\overline De^{-t\overline D^2}$ on
  $\overline M$ using the fixed spin structure and the metric
  $\Psi^*g$, whereas $\kappa_t(x,y)$ is the same function on
  $\Psi^*\overline M$ defined using the pullback structure throughout.

  In particular,  reasoning as in the proof of
 Propoistion  \ref{prop:diffeo_and_rplus},
 we see that 
  for $h\in\Gamma$
  \begin{equation}\label{eq:diffeo_formula}
    \eta_{<h>}(\overline D_{\Psi^*g}) = \eta_{<\alpha_\Psi(h)>}(\overline D),
  \end{equation}
  whenever $\eta_{<h>}$ is defined.
%
\end{remark}

The following example  is a direct consequence of the results of Botvinnik and Gilkey \cite{Gilkey}.
\begin{example}\label{ex:non_trivial_rho_for_Zn}
  Let $\integers/n$ be a finite cyclic group, and $m>4$ be congruent $3$
  mod $4$. Then $\rho_{(2)}\colon \Pos_m^{\spin}(B\integers/n)\to\reals$
  is non-trivial.
  Since it is a group homomorphism for the additive group of $\reals$,
  its image is infinite. 
\end{example}
\begin{proof}
  We only have to observe that $\rho_{(2)}$ is a
  twisted rho-invariant, where we twist with
  $-\reals+\frac{1}{n}\reals[\integers/n]$. Indeed, the first representation is
  the opposite of the trivial representation, giving minus the ordinary
  eta-invariant; the second one is a multiple of the regular representation,
  giving the $L^2$-eta invariant. In order to prove the last statement recall that
   for any unitary representation $\phi$ with character $\chi_\phi$,
  the twisted eta invariant $\eta_\phi (D)$ can be expressed  by
  \begin{equation}\label{fourier}
  \eta_\phi (D)=
  \sum_{h\in\integers/n} \chi_\phi (h) \eta_{h} (\overline{D})\,
  \end{equation}
  where $<h>=h$ given that the group is commutative.
  Since the character of the regular representation is the delta function at
  the identity element,
  we see that the eta invariant for the operator twisted by the regular representation
  is nothing but  the $\eta$-invariant of the $\integers/n$-covering, which is
  $n$-times the $L^2$-eta invariant of this covering.\\ The
  character $\chi$ of the virtual representation $-\reals+1/n
  \reals[\integers/n]$ is invariant under
  inversion: $\chi(g)=\chi(g^{-1})$. This means, by definition,  that  
$-\reals+\frac{1}{n}\reals[\integers/n]\in R^+_0 (\integers/n)$ where
\begin{equation}\label{def-of-r+}
R^+_0 (\integers/n):= \{\phi\in R(\integers/n)\,|\,\chi_\phi(1)=0\,; \chi_\phi(g)=
  \chi_\phi (g^{-1})\,\,\forall g\}.
  \end{equation}
  By the results of  Botvinnik-Gilkey \cite[Proof of Theorem 2.1]{BoGi}
  we know that 
  \begin{equation}\label{bo-gi-precise}
  \forall \phi\in R^+_0 (\integers/n) \,\,\exists\,\, [M,g,u] \in \Pos_m^{\spin}(B\integers/n)\,\,
\text{such that}\,\, \rho_\phi [M,g,u]\not= 0
\end{equation}
and it suffices to apply this result to $\reals-\frac{1}{n}\reals[\integers/n]$ .
\end{proof}


\begin{remark}\label{remark:bot-gil}
Let $\Gamma$ be any  {\bf finite}  group and let $m>4$ be congruent 3 mod 4.
The results of Botvinnik and Gilkey show  that the map
$\Psi: R^+_0 (\Gamma)\otimes \CC \rightarrow ({\rm Pos}^{\spin}_m(B\Gamma)\otimes \CC)^\prime$
defined by
\begin{equation}\label{dual-valued}
\Psi(\phi)[M,g,u]:=\rho_\phi [M,g,u] 
\end{equation}
is injective.
Let ${\rm Class}(\Gamma)=\{f:\Gamma\to\CC\,\,|\,\,f(\gamma^{-1} h \gamma)=f(h)\,\,\forall \gamma, h\in \Gamma\}$ be the complex vector space
of class functions on $\Gamma$. Let 
\begin{equation}\label{class+}
{\rm Class}^+_0 (\Gamma)=\{f\in {\rm Class}(\Gamma)\,\,|\,\,
f(1)=0, f(h)=f(h^{-1})\}\,.
\end{equation}
 Then there is a natural isomorphism of vector spaces
$\Theta: R^+_0 (\Gamma)\otimes\CC \to {\rm Class}^+_0 (\Gamma)$ obtained by associating to 
$\phi\in  R^+_0 (\Gamma)$ its character $\chi_\phi$.
There is also a map 
$\Phi: {\rm Class}^+_0 (\Gamma) \rightarrow ({\rm Pos}^{\spin}_m (B\Gamma)\otimes \CC)^\prime$
given by 
\begin{equation}\label{dual-valued-2}
\Phi (f) [M,g,u]:= \sum_{<h>} \rho_{<h>} [M,g,u] f(<h>)\,.
\end{equation}
Since by the analog  of $\eqref{fourier}$ we see that $\Phi\circ \Theta = \Psi$,
we conclude that $\Phi$ is also injective if $m>4$ is congruent 3 mod 4.\\
We shall apply this result to $\Gamma=\ZZ/n$: thus for these values of $m$
\begin{equation}\label{non-vanish-cyclic}
\forall f\in {\rm Class}^+_0 (\ZZ/n)\, \exists\, y \in {\rm Pos}^{\spin}_m (B\ZZ/n)\otimes \CC\,:\quad\sum_{h} \rho_{h} (y) f(h)\not=0
\end{equation}
\end{remark}

The following lemma describes how to compute delocalized rho-invariants for
manifolds obtained by induction.

\begin{lemma}\label{lem:rho_2_and_induction}
  Let $\pi=\integers/n$ be a finite cyclic group, $j\colon \integers/n\into
  \Gamma$ 
  an injective group homomorphism. Fix $1\ne g\in \Gamma$ of finite order and
  such that its conjugacy class $<g>$ has polynomial growth. The
  delocalized rho-invariant $\rho_{<g>}\colon \Pos^{\spin}_m(B\Gamma)\to
  \reals$ is defined if $m$ is odd. We have the induced map $Bj_*\colon
  \Pos^{\spin}_m(B\integers/n)\to \Pos^{\spin}_m(B\Gamma)$. Then
  \begin{equation}\label{eq:induction}
    \rho_{<g>}(Bj_* x) = \sum_{h\in j^{-1}(<g>)} \rho_{<h>}(x);\qquad\forall
    x\in \Pos^{\spin}_m(B\integers/n).
  \end{equation}

  Similarly, considering the $L^2$-rho-invariant
  \begin{equation*}
    \rho_{(2)}(Bj_* x) =   \rho_{(2)}(x);\qquad\forall
    x\in \Pos^{\spin}_m(B\integers/n).
  \end{equation*}
\end{lemma}
\begin{proof}

  This is a well known feature of $L^2$-invariants. We
  indicate the proof, showing along the way how it extends to the delocalized
  invariants. 
  Assume that $x=[M,g,u\colon M\to B\pi]$. Observe that $j$ is
  injective. This implies that the covering $(Bj)^*E\Gamma\to
  B\pi$ decomposes as a 
  disjoint union (parametrized by the elements of the set
  $\Gamma/j(\pi)$) of copies of $E\pi$. For the convenience of the
  reader we recall a possible argument. Given the universal free
  $\Gamma$-space $E\Gamma$, the action
  of $\pi$ on $E\Gamma$ (via $j$) allows us to view $E\Gamma$ as a
  model of $E\pi$, with $B\pi:=E\Gamma/\pi$. In this picture, $Bj$ is
  simply the projection map $E\Gamma/\pi\to E\Gamma/\Gamma$. Then the
  pullback $(Bj)^*E\Gamma=\{(x\pi,x\gamma)\in E\Gamma/\pi\times E\Gamma\mid x\pi\in
  E\Gamma/\pi,\gamma\in\Gamma\} \iso E\Gamma\times \pi\backslash\Gamma$
  with the evident map $(x\pi,x\gamma)\mapsto (x\pi,\pi\gamma)$.

Consequently, the covering $\overline M=(Bj\circ
  u)^* E\Gamma=u^*(Bj)^*E\Gamma$ decomposes as a disjoint union of copies of the covering
  $\tilde M$ classified by $u$. The construction of the
  $L^2$-eta invariant for this disjoint union $\overline M=(Bj 
\circ u)^* E\Gamma$ involves only the one component $\tilde M$
which contains the fundamental domain, and therefore is
exactly the same as the construction of the $L^2$-eta invariant for
$\tilde M$ itself. Since the ordinary $\eta$-invariant does only depend on
  $M$, also the $L^2$-rho invariants coincide. 

More precisely, and also
  holding for the delocalized 
  invariants, recall from Definition \ref{def:rho_2} that
  \begin{equation*}
    \rho_{<h>}(\overline D) = \frac{1}{\sqrt{\pi}}\int_{0}^\infty
    \sum_{\gamma\in<h>} \int_{\mathcal{F}} \tr_x k_t(x,\gamma
    x)\,dx\,\frac{dt}{\sqrt{t}}.
  \end{equation*}
  Now $\overline M$ decomposes as a disjoint union of copies of $\tilde
  M$. The heat kernel $k_t(x,y)$ 
  vanishes if $x$ and $y$ belong to different components, and if $x$ and $y$
  lie in the same component, coincides with
  the heat kernel of the operator restricted to that component (use uniqueness
  of the heat kernel).
  If $x\in \mathcal{F}\subset \tilde M$ and $\gamma\in \im(j)$ then $\gamma
  x\in \tilde M$ (because $\tilde M\to M$ is just the covering corresponding
  to the subgroup $j(\pi)$ of $\Gamma$). However, if $\gamma\notin \im(j)$,
  then $\gamma x\notin \tilde M$ (for the same reason). Thus, in the sum
  above, all summands with $\gamma\notin \im(j)$ vanish, whereas the summands
  with $\gamma\in\im(j)$ are exactly those (using an obvious diffeomorphism)
  showing up in the definition of the delocalized invariants for $\tilde D$ on
  $\tilde M$, and this is what is stated in the assertion of the Lemma.
\end{proof}

\begin{remark}
  The proof of Lemma \ref{lem:rho_2_and_induction} gives also a
  formula for induction from arbitrary (not necessarily cyclic)
  subgroups. Namely, if $j\colon \pi\into \Gamma$ is an injective homomorphism
  for a not necessarily finite cyclic group,
  \begin{equation*}
    \rho_{<h>}(Bj_*(x)) = \sum_{<\gamma>\subset j^{-1}(<h>)}
    \rho_{<\gamma>}(x), 
    \end{equation*}
    where the sum on the right hand side runs over all the $\pi$-conjugacy
    classes which are contained in $j^{-1}(<h>)$.
\end{remark}

\begin{theorem}\label{theo:infinitely-many}
  Assume that $M$ is a spin-manifold of dimension $m>4$, $m\congruent 3
  \pmod 4$. Assume that $g$ is a metric with positive scalar curvature
  on $M$, and that the fundamental group $\Gamma$ of $M$ contains at
  least one non-trivial element of finite order.  Then $\charPos(M,u)$
  is infinite, i.e.~$M$ admits infinitely
  many different \emph{bordism} classes of metric with $\scal>0$. They are
  distinguished by $\rho_{(2)}$.

  More precisely, the infinitely many bordism classes we construct are also
  different modulo the ``action'' of the diffeomorphism group, i.e.~we get
  metrics $(g_\alpha)_{\alpha\in A}$ such that $\abs{A}=\infty$ and for every
  diffeomorphism $f$ of $M$,
  $f^*g_\alpha$ is bordant to $g_\beta$ only if $\alpha=\beta$.

As a consequence, the space $\mathcal{R}^+(M)/\Diffeo(M)$, the moduli space of
  metrics of positive scalar curvature, has infinitely many components,
  distinguished by $\rho_{(2)}$.
\end{theorem}
\begin{remark}
  Recall that by the methods of Gromov and Lawson, it is known that
  $\mathcal{R}^+(M)/\Diffeo(M)$ has infinitely many components for every
  manifold of dimension $4k+3$, $k\ge 1$ (compare
  \cite[Theorem 7.7]{lawson89:_spin}). 
 Strictly speaking the result stated in  \cite{lawson89:_spin}
  only involves   $\mathcal{R}^+(M)$:
  an inspection of the proof shows that  the main argument used there
  also establishes the fact that $|\pi_0 (\mathcal{R}^+(M)/\Diffeo(M))|$ $=$$\infty$
  : indeed it suffices to observe that the signature is a cut-and-paste invariant.
   Notice  however, that
 by construction the examples they get are
  all bordant to each other.
\end{remark}
\begin{proof} Let $j\colon \integers/n\to \Gamma$ be an
  injection. This exists for some $n>1$ since $\Gamma$ is not torsion free.
  By Example \ref{ex:non_trivial_rho_for_Zn}, the homomorphism $\rho_{(2)}\colon
  \Pos_m^{\spin}(B\integers/n)\to\reals$ is non-trivial (therefore
  has infinite image). The group
  $\Omega_m^{spin}(B\integers/n)$ is finite by the Atiyah-Hirzebruch
  spectral sequence. Consequently, the kernel $K$ of the map
  $\Pos_m^{\spin}(B\integers/n)\to \Omega^{\spin}_m(B\integers/n)$ has
  finite index, and the restriction $\rho_{(2)}|\colon K\to \reals$
  also is non-trivial with infinite image.

  Let $u\colon M\to B\Gamma$ be the chosen classifying map of a universal
  covering. By naturality of the exact sequence
  \eqref{eq:pos_sequemnce} and Proposition
  \ref{prop:action_of_ker_on_charPos} $Bj_* k+[M,g,u] \in \charPos(M,u)$
  for each $k\in K$. Moreover, by Lemma \ref{lem:rho_2_and_induction},
  \begin{equation*}
\rho_{(2)}(Bj_*k +[M,g,u]) =
  \rho_{(2)}(k) + \rho_{(2)}(M,g,u),
\end{equation*}
  Consequently, $\rho_{(2)}\colon \charPos(M,u)\to\reals$ has infinite
  image.

  Using Proposition \ref{prop:diffeo_and_rplus} and the surjectivity of
  $\mathcal{R}^+(M)\to \charPos(M,u)$, the map $\rho_{(2)}\colon
  \pi_0(\mathcal{R}^+(M)/{\rm Diffeo}_{\#}(M))\to \reals$ also has infinite image.

  Since the spin-structure preserving
  diffeomorphisms have finite index in all diffeomorphisms, even
  modulo $\Diffeo(M)$ there are infinitely many components in the moduli
  space. In a similar way, we can get infinitely many bordism classes which
  are different even modulo pullback with arbitrary diffeomorphisms.
\end{proof}

\subsection{Different conjugacy classes of torsion elements in the fundamental
  group and positive scalar curvature}
\label{sec:diff-conj-class}

In this subsection we shall sharpen Theorem \ref{theo:infinitely-many}
and extend it to dimensions $4k+1$
under some additional assumptions on $\Gamma$.
 
\begin{theorem}\label{theo:infinitely-many-sharp}
  Let $\Gamma$ be a discrete group. Consider the following subset of 
  the set $\C$ of all conjugacy classes of $\Gamma$:
  \begin{equation*}
\C_{fp}:=\{<h>\subset \Gamma\mid h\text{ has finite order }, <h>\text{ has
  polynomial growth}\}.
\end{equation*}
On this set, we have an involution $\tau$ given by $<h>\mapsto <h^{-1}>$.
Assume that $M$ is a closed spin manifold with fundamental group $\Gamma$,
with classifying map $u\colon M\to B\Gamma$, of dimension $4k+3$, $k\ge 1$.
Then, on $\charPos(M,u)$, a free abelian group of rank $\abs{\C_{fp}/\tau}$ acts
freely.
\end{theorem}
\begin{proof}
Let $K:=\ker(\Pos^{\spin}_m(B\Gamma)\to
  \Omega_m^{\spin}(B\Gamma))$: then
it suffices to show that 
$\dim K\otimes \CC\geq \abs{\C_{fp}/\tau}$. 
For each $<g>\in \C_{fp}$ consider the characteristic function $\kappa(g)$
of the set $<g>\cup <g^{-1}>$. By \eqref{class+}, $\kappa(g)$ belongs to ${\rm Class}^+_0 (\Gamma)$.
Let $L_{fp}$ be the vector subspace of  ${\rm Class}^+_0 (\Gamma)$ whose elements are finite
linear combinations
of  $\kappa(g)$, with $<g>\in \C_{fp}$. This is a vector space of dimension 
$\abs{\C_{fp}/\tau}$ and we denote by $\kappa$, $\kappa=\sum_j \lambda_j \kappa(g_j)$,
the generic element.
Following Remark \ref{remark:bot-gil}, we begin by showing that the map
$\Phi: L_{fp}\subset {\rm Class}^+_0 (\Gamma)\rightarrow ({\rm Pos}^{\spin}_m (B\Gamma)\otimes \CC)^\prime$
which associates to $\kappa=\sum_j \lambda_j \kappa(g_j)$ the functional $\Phi(\kappa)$,
\begin{equation}\label{eq:def_of_Phi}
  \Phi (\kappa)[M,g,u]:=\sum_j \lambda_j (\rho_{<g_j>}[M,g,u]  + \rho_{<g_j^{-1}>}[M,g,u])\,,
\end{equation}
is injective.
Choose $g_\ell$ so that $\lambda_\ell\not= 0$. Let $\pi(g_\ell)$ be the finite cyclic group
generated by $g_\ell$. Consider the restriction $\kappa |_{\pi(g_\ell)}$, 
an element in $ {\rm Class}^+_0 (\pi(g_\ell))$. Then by the results of Botvinnik-Gilkey,
as stated in \eqref{non-vanish-cyclic}, we know that there exist $y\in {\rm Pos}^{\spin}_m (B\pi(g_\ell))\otimes \CC$ such that 
\begin{equation*}
\sum_{h\in \pi(g_\ell)}\rho_{h} (y) \kappa |_{\pi(g_\ell)} (h)\not=0.
\end{equation*}
Let $j: \pi(g_\ell)\hookrightarrow \Gamma$ be the natural inclusion and let $x:=Bj_* (y)$
so that $x\in {\rm Pos}^{\spin}_m (B\Gamma)\otimes \CC$. By the
induction formula \eqref{eq:induction} we know   that
\begin{equation*}
  \Phi (\kappa)(x)
  = \sum_{h\in \pi(g_\ell)}\rho_{h} (y) \kappa |_{\pi(g_\ell)} (h)
  \end{equation*}
and we can therefore conclude that $ \Phi (\kappa)(x)\not=0$.
It remains to show that we can choose $x\in \ker(\Pos^{\spin}_m(B\Gamma)\to
  \Omega_m^{\spin}(B\Gamma))\otimes \CC$. By naturality of the long exact sequence
  \eqref{eq:pos_sequemnce} it suffices to show that we can choose
  $y\in \ker (\Pos^{\spin}_m(B\pi(g_\ell))\to
  \Omega_m^{\spin}(B\pi(g_\ell)))\otimes \CC$. However, since
  $ \Omega_m^{\spin}(B\pi(g_\ell))$ is finite, this is easily accomplished by
  taking a suitable multiple of the original $y$.
  \end{proof}






We now analyze the case $\dim(M)\equiv 1\pmod 4$.
Let
\begin{align}\label{def-of-r_-and-class-}
&R^-_0 (\integers/n):= \{\phi\in R(\integers/n)\,|\,\chi_\phi(1)=0\,; \chi_\phi(h)=-
  \chi_\phi (h^{-1})\,\,\forall h\}\\
& {\rm Class}^-_0 (\Gamma)=\{f\in {\rm Class}(\Gamma)\,\,|\,\,
f(1)=0, f(h)=-f(h^{-1})\}\,.
\end{align}
Then the results of Botvinnik and Gilkey in \cite{BoGi}
imply that the analogs of 
\eqref{bo-gi-precise} of Example \ref{ex:non_trivial_rho_for_Zn} 
and of \eqref{non-vanish-cyclic} of Remark \ref{remark:bot-gil}
hold.
For the convenience of the reader we explicitly restate the latter property 
in this
new context:
\begin{equation}\label{non-vanish-cyclic-minus}
\forall f\in {\rm Class}^-_0 (\ZZ/n)\, \exists\, y \in {\rm Pos}^{\spin}_m (B\ZZ/n)\otimes \CC\,\,
|\,\,\sum_{h} \rho_{h} (y) f(h)\not=0
\end{equation}

\begin{theorem}\label{theo:diff-conj-class-1}
  Let $\Gamma$ be a discrete group.
Let   $\C_{fp}$ and $\tau: \C_{fp}\to \C_{fp}$ be as in the statement of Theorem
 \ref{theo:infinitely-many-sharp}. \\ Let
 $\C^0_{fp}=\{<h>\in \C_{fp} \,|\, <h>\not= <h^{-1}>\}$. 
  Assume that $M$ is a closed spin manifold with fundamental group $\Gamma$,
with classifying map $u\colon M\to B\Gamma$, of dimension $4k+1$, $k\ge 1$.
Then, on $\charPos(M,u)$, a free abelian group of rank $\abs{\C^0_{fp}/\tau}$ acts
freely.

Moreover, if $\C_{fp}^0$ is not empty then there are infinitely many bordism classes which are
different modulo the ``action'' of the diffeomorphism group as in
Theorem \ref{theo:infinitely-many}.
As a consequence, the space $\mathcal{R}^+(M)/\Diffeo(M)$, the moduli space of
  metrics of positive scalar curvature, has infinitely many components
  in our situation,
  distinguished by the collection $\rho_{<h>}$, $<h>\in \C^0_{fp}$.
\end{theorem}

\begin{proof}
Let $<h>\in \C^0_{fp}$ and consider the function $\kappa(h)$
which is equal to $1$ on $<h>$, equal to $-1$ on $<h^{-1}>$
and $0$ elsewhere; $\kappa(h)$ so defined is an element of
${\rm Class}^-_0 (B\Gamma)$.  Let $L^0$ be the vector subspace 
of ${\rm Class}^-_0 (B\Gamma)$ whose elements are finite linear
combinations of $\kappa(h)$, with $<h>\in \C^0_{fp}$. This is a vector space
of dimension $|\C^0_{fp}/\tau|$.
Using the induction formula and \eqref{non-vanish-cyclic-minus} the proof 
now proceeds as in the proof of Theorem \ref{theo:infinitely-many-sharp}.

If $\C^0_{fp}$ is not empty,
choose the collection of functions $\kappa(h)$ of $L^0$ for $<h>\in
\C^0_{fp}$. Note that $\kappa(h^{-1})=-\kappa(h)$. Then, dualizing $\Phi$, we get a
 map
\begin{equation*}
 \mathcal{R}^+(M)\to  \textrm{Pos}^{\spin}_m (B\Gamma)\otimes \CC \to \complexs^{\C^0_{fp}};
[M,g,u]\mapsto (\Phi(\kappa(h))(M,g,u))_{<h>\in \C^0_{fp}}
\end{equation*}
with infinite image. If we had chosen one half of the functions
$\kappa(h)$, forming a basis, the map would have been surjective.

Now, given a spin structure preserving diffeomorphism $\Psi\colon M\to
M$ (with a given lift to the spin principal bundle), we get an induced
automorphism $\alpha_\Psi$ of $\Gamma$ as in the proof of Proposition
\ref{prop:diffeo_and_rplus}, and an
induced permutation of $\C^0_{fp}$. Moreover, by
\eqref{eq:diffeo_formula},
\begin{equation*}
  \Phi(\kappa(h))(M,\Psi^*g,u) = \Phi(\kappa(\alpha_\Psi(h)))(M,g,u)
\end{equation*}
so that we above map induces a well defined map
\begin{equation*}
  \pi_0(\mathcal{R}^+(M)/\Diffeo_\#(M)) \to \CC^{\C^0_{fp}}/\Sigma,
\end{equation*}
where we quotient the right hand side by the action of the permutation
group, permuting the entries of the vector. Since this group is
finite, the image still is infinite.

Since the spin structure preserving diffeomorphisms have finite
index in the whole diffeomorphism group, even
$\pi_0(\mathcal{R}^+(M)/\Diffeo(M))$ is infinite.
\end{proof}

\begin{remark}
Notice, in particular, that if $\dim(M)\equiv 1\pmod 4$ and if $\Gamma$ contains an element $g$ of finite order not conjugate to
  its inverse and such that the conjugacy class $<g>$ has polynomial
  growth, then a manifold $M$ as above admits infinitely many pairwise non-bordant metrics
  of positive scalar curvature. To our knowledge, this is the first such
  result of considerable generality.
\end{remark}

\begin{remark}
  We want to point out that there are many non-trivial examples of groups $\Gamma$, where
  $\mathcal{C}^0_{fp}$ is non-empty. In particular, this applies to
  \begin{enumerate}
  \item Groups with a central element of odd order (here the relevant
  conjugacy class consists of one element). For an arbitrary group $H$ and a
  finite group $F$ (of odd order), all non-trivial elements of $F$ in $F\times H$ have
  this property. 
  \item Many groups with a non-trivial finite conjugacy center, consisting of
  elements of finite order. Such groups are e.g.~obtained as extensions $1\to
  F\to G\to H\to 1$ with $F$ finite (and of odd order).
\item groups of polynomial growth with elements of finite order (in this case,
  every conjugacy class has of course polynomial growth).
\item The restricted wreath product $(\oplus_{k\in\integers}\integers/n)\semiProd \integers$ (if $n=2$
  this is called the lamplighter group) is a group of exponential growth, such
  that every element in the normal subgroup
  $\oplus_{k\in\integers}\integers/n$ has an infinite conjugacy class of
  polynomial growth.
\end{enumerate}
Similar examples give rise to non-empty $\mathcal{C}_{fp}$.
\end{remark}

\medskip

\noindent
{\bf Further questions and open problems}:

\smallskip

1) We study only $\pi_1$-bordism, which is necessary for our method, 
because it
uses the common fundamental group throughout. Nonetheless, this concept is
somewhat unnatural from a geometric point of view. It would be interesting 
to
know whether our examples remain non-bordant if we talk about the most 
obvious
simple definition of bordism of metrics of positive
scalar curvature, or to find any examples which are not bordant in this 
weak
sense.

2) We get some information about the number of components of the space of
   metrics of positive scalar curvature. What else can be said about its
   topology, in particular about higher homotopy groups?

3) We prove that for spin manifolds of dimension $4k+1$ with positive scalar
   curvature and with fundamental group which contains a central element 
of   odd order, the moduli space of metrics of positive scalar curvature has
   infinitely many components. In dimension $4k+3$ this is known 
unconditionally
   - what about the given dimension $4k+1$.

\section{An example of a non-convergent delocalized eta invariant}
\label{sec:an-example-non}

 In this  section we  compute Lott's delocalized $\eta$-invariant of an easy example, and use
  this to produce an example where it does not converge.

\medskip
  Consider the manifold $S^1$ with the usual metric. The Dirac (and signature)
  operator of $S^1$ is (unitarily equivalent to) the operator $D=\frac{1}{i}\frac{d}{dx}$.

  The integral kernel $k_t(x,y)$ of $\tilde D\exp(-t\tilde D^2)$ on the
  universal 
  covering $\reals$ of $S^1$ is
  \begin{equation*}
    k_t(x,y)=-i\frac{x-y}{2t\sqrt{2\pi t}} e^{\frac{(x-y)^2}{4t}}.
  \end{equation*}
  Fixing the fundamental domain $\mathcal{F}=[0,1]$ for the covering
  projection, and using the action by the deck transformation group
  $\integers$ by addition: $(x,n)\mapsto x+n$,
  the delocalized eta-invariant for a subset $X\subset \naturals$ formally
  would be 
  \begin{equation*}
    \begin{split}
      \eta_X(\tilde D):=& \frac{1}{\sqrt{\pi}}\int_0^\infty \sum_{n\in X}
      \int_{\mathcal{F}} k_t(x,x+n)\, dx\,\frac{dt}{\sqrt{t}}\\
      =&
      -\frac{i}{4{\pi}}\int_0^\infty \sum_{n\in X} \int_0^1 \frac{n}{t^2}
      e^{n^2/4t}\,dx\,dt\\
      =& -\frac{i}{4\pi} \sum_{n\in X}\frac{1}{n}  \int_0^\infty \frac{e^{1/4t}}{t^2}\,dt,
  \end{split}
\end{equation*}
where at the end we use the substitution $t/n^2=s$ and the fact that the
integrands are all positive, such that we can interchange the summation over
$n\in X\subset \naturals$ and the integral over $t$.

It is clear that this expression is divergent for suitable infinite $X\subset \naturals$.

Consider next the group $\Gamma=\rationals\semiProd \left(\bigoplus_{n\in\integers}\integers\right)$, where
the generator of the $n$-th summand of $\bigoplus\integers$ acts by
multiplication with the $\abs{n}$-th prime number. By
the definition of semidirect products, the conjugacy class of $1\in\rationals$
in the 
kernel group is exactly $\rationals_{>0}$. Its intersection with the
subgroup $\integers$ generated by $1$ is therefore $
\naturals_{>0}\subset\integers$.

Consider also $G:= \Gamma \semiprod_{\alpha}$, the HNN-extension of $\Gamma$
along  $\alpha\colon \bigoplus_{n\in\integers}\integers\to
\bigoplus_{n\in\integers}\integers;( n\mapsto \lambda_n)\mapsto (n\mapsto \lambda_{n+1})$,  the shift of the non-normal subgroup
$\bigoplus_{n\in\integers}\integers$. Then $G$ is generated by $3$ elements:
 $1$ in the additive groups of $\rationals$, a generator of the copy of
 $\integers$ labelled with zero in $\bigoplus_{n\in\integers}\integers$, and
 the stable letter $t$. Moreover, using the normal form of elements in an
 HNN-extension, the intersection of the conjugacy class of $1$ with
 $\rationals$ still constists of $\rationals_{>0}$, and therefore the
 interesection with the additive subgroup of integers consists of the natural
 numbers.

Observe that $G$ is finitely generated, but by its definition only recursively
countably presented. As such, by a standard procedure, $G$ can be embedded
into a finitely presented group $H$ wich is obtained as follows (compare
\cite[Theorem 12.18]{MR1307623}).

One first constructs an auxiliary group $B_2$, then considers the group
$B_3=B_2*G$, the free product of $B_2$ and $G$. The next group is an
HNN-extension of $B_3$ along a subgroup which is of the form $U*G$ for a
suitable subgroup $U$ of $B_2$. The stable letters act trivially on $G$. By
the normal form of elements of an HNN-extension, it follows that for every
element $x\in G$, the conjugacy class of $x$ in $G$ is equal to the
intersection of the conjugacy class of $x$ in $B_3$ with $G$.

In the next steps, one constructs two further HNN-extensions of the previous
group (starting with $B_3$) along subgroups with trivial intersection with
$G$. Again, it follows from the normal form of elements in an HNN-extension
that for every $x\in G$ the conjugacy class of $x$ in $G$ coincides with the
intersection of $G$ with the conjugacy class of $x$ in the bigger group. The
final group $H:=B_6$ is finitely presented, contains $G$ (and therefore
$\integers$) as a subgroup, and the intersection of the conjugacy class of $1$
with $\integers$ consists exactly of the positive integers.

Consider $u\colon S^1\to B\integers\to B\Gamma\to BG\to BH$, where the first map is the
classifying map for the universal covering (i.e.~the identity if we use the
model $B\integers=S^1$,) and the other maps are induced by the inclusion
$\integers\into \Gamma\into G\into H$ (the first inclusion sends $1\in\integers$ to $1\in\rationals\subset\Gamma$).

Let $\overline M\to S^1$ be the induced covering, and $\overline D$ the lift
of $D$ to this covering. Then, by the formula for delocalized eta-invariants of
induced manifolds,
\begin{equation*}
  \eta_{<1>}(\overline D) = \eta_{P}(\tilde D)
\end{equation*}
with $P=\naturals_{>0}\subset\integers$, which is not convergent by the above
calculation.

This is an example of an operator where the delocalized eta-invariant of John
Lott is not defined.

\begin{remark}
  The same calculation works for the product of a manifold $M$ of
  dimension   $4k$ with $S^1$ with product metric. During the calculations,
  one has to multiply the above expressions for $S^1$ with $\hat{A}(M)$. If
  this number is non-zero, one therefore gets the same non-convergence
  behaviour for manifolds of arbitrarily high dimension.

 Similar calculations should also be possible for  more general mapping tori of
   a $4k$-manifold, compare \cite{MR2000k:58039}. One should be able to work
  with the signature as well as the Dirac operator.
\end{remark}

\begin{remark}
  It is probably not trivial to obtain an example where the conjugacy-class
  (inside the new group $\Gamma$) has polynomial growth. Observe that this is
  not the case for the construction we describe.

  It would also be very interesting to find examples of non-convergence with
  positive scalar curvature (then, necessarily, the conjugacy class could not have
  to have polynomial growth).

  It would be even more interesting if one could produce examples as above
  where the fundamental group of the manifold is the group $H$. It is not
  clear to us how to construct such an example and keep control of the
  calculation of the $\eta$-invariant.

  Another open problem is the construction of examples with non-convergent
  delocalized $L^2$-Betti numbers. As a starting point, one should again look
  for manifolds with many non-trivial such; by induction to larger groups one
  might then be able to obtain one conjugacy class where the invariants don't
  converge. 
\end{remark}

 \section{Possible values of APS-rho invariants for the Dirac operator}
\label{sec:possible-values-aps}

In this section, we prove Proposition
\ref{prop:possible_values_of_APS_rho_in_general}. Its proof is
modeled on the proof of the corresponding statement \cite[Theorem
7.1]{Hi-Gu-We} for the signature operator.

Let $M$ be a closed spin manifold with positive scalar curvature. Let
$u\colon M\to B\Gamma$ be a continuous map and
$\lambda_1,\lambda_2\colon \Gamma\to U(d)$ two finite dimensional
unitary representations of $\Gamma$. Set $\Gamma_1:=\im(\lambda_1)$
and $\Gamma_2:=\im(\lambda_2)$. We consider $\Gamma_1$ and $\Gamma_2$
as \emph{discrete} groups which happen to be subgroups of $U(d)$.

We compose $u$ with the maps induced by $\lambda_{1}$ and $\lambda_2$
to get $v\colon M\to B[\Gamma_1\times\Gamma_2]$. The tuple
$(M, v)$ then represents an element $[M,v]$ in the real K-homology of
$B[\Gamma_1\times\Gamma_2]$. We can now apply the reduced Baum-Connes
map $\mu_{red}$ to this element, to get
\begin{equation*}
  \ind(D_L) \in KO_*(C^*_{red}(\Gamma_1\times\Gamma_2)).
\end{equation*}
Here $L$ is the Mishchenko-Fomenko line bundle associated to $v$.

Since $M$ has positive scalar curvature, this index is zero by the
Lichnerowicz formula. On the other hand, $\Gamma_1\times \Gamma_2$ is
a linear group by its very construction. By the main result (0.1) of
\cite{Hi-Gu-We}, the following map
\begin{equation*}
  \mu_{red}\colon
  K^{\Gamma_1\times\Gamma_2}_*(E[\Gamma_1\times\Gamma_2]) \to
  K_*(C^*_{red}(\Gamma_1\times\Gamma_2))
\end{equation*}
is split injective in this case. Their proof  applies in the same way
to the real K-theory, since they really prove that linear groups
uniformly embed into Hilbert space, which implies the coarse
Baum-Connes isomorphism conjecture for linear groups. This in turn implies the
real coarse Baum-Connes conjecture for linear groups by a well known
principle, compare e.g.~\cite{math.KT/0311295}. From here, the descent principle
implies split injectivity of the usual real reduced Baum-Connes
map. Therefore
\begin{equation*}
  KO^{\Gamma_1\times\Gamma_2}_*(E[\Gamma_1\times\Gamma_2]) \to
  KO_*(C^*_{\reals,red}(\Gamma_1\times\Gamma_2))
\end{equation*}
is also injective.

We have to produce a link between $
KO^{\Gamma_1\times\Gamma_2}_*(E[\Gamma_1\times\Gamma_2])$ and  the
non-equivariant Baum-Connes map with $KO_*(B[\Gamma_1\times\Gamma_2])$
used
so far. There is a canonical map
\begin{equation*}
  KO_*(B[\Gamma_1\times\Gamma_2]) \to
  KO_*^{\Gamma_1\times\Gamma_2}(E[\Gamma_1\times\Gamma_2])
\end{equation*}
By standard arguments (compare \cite[Section 7]{Hi-Gu-We} and
\cite[Lemma 2.9]{Lueck-Stamm}) in equivariant homology theory, this map is
split injective after tensoring with $R$ of
Proposition \ref{prop:possible_values_of_APS_rho_in_general}.

Putting these two facts together,
\begin{equation*}
  \mu_{red}\colon KO_*(B[\Gamma_1\times\Gamma_2])\tensor R\to
  KO^*(C^*_{red}(\Gamma_1\times\Gamma_2))\tensor R
\end{equation*}
is injective. Since we have already seen that $\ind(D_L)=0$, this
implies that there is $l\in \naturals$ which is a product of orders of
elements of $\Gamma_1\times\Gamma_2$ such that
\begin{equation*}
 l\cdot [M,u]=0 \in KO_*(B[\Gamma_1\times\Gamma_2]).
\end{equation*}

We now use the geometric
description of $KO_*(X)$ in terms of spin bordism due to Hopkins-Hovey \cite[Theorem 1]{Ho-Ho}. 
First observe that there is a natural map 
$\Omega^{spin}_*(X)
\rightarrow KO_*(X)$
which assigns to a spin manifold $M$ with map $v\colon M\to X$ the
class $[M,v]\in KO_*(X)$ given by the geometric description of
$KO_*(X)$.
  Next, consider the special case  $\tau\colon\Omega^{spin}_*(pt)\to KO_*(pt)$
  of this 
  homomorphism for $X$ equal to a point. This is a  (graded)
  ring homomorphism with  kernel consisting  of  some 
  manifolds with vanishing $\hat{A}$-genus, and  cokernel   the
  ideal generated by 
  $KO_{-8}(pt)$. We can consider $KO_*(pt)$ as a module
  over $\Omega^{spin}_*(pt)$ via  $\tau$ 
 and form
  $\Omega^{spin}_*(X)\tensor_{\Omega^{spin}_*(pt)} KO_*(pt)$. The result of Hopkins and Hovey
  says that for each CW-complex $X$ the induced map is an isomorphism: 
\begin{equation*}
\Omega^{spin}_*(X)\tensor_{\Omega^{spin}_*(pt)} KO_*(pt)
\xrightarrow{\iso} KO_*(X).
\end{equation*}

\begin{definition}
 A  Bott manifold, $B$ is a $8$-dimensional simply connected spin manifold
  with $\hat{A}(B)=1$.
\end{definition}

\begin{lemma}
  If $[M,v]=0\in KO_*(X)$ then there are $n\in\naturals$ and spin
  manifolds $A_i,C_i$ with $\hat{A}(C_i)=0$ and continuous maps
  $u_i\colon A_i \to X$ such that $[M,v]\times B^n$ is bordant in
  $\Omega^{spin}_*(X)$ to the disjoint union of $[A_i\times
  C_i,u_i]$. The maps to $X$ are given by first projecting to the
  first factor and then using $v$ or $u_i$, respectively.
\end{lemma}
\begin{proof}
  Considering  $KO_*(pt)$ as a module
  over $\Omega^{spin}_*(pt)$ as above, we obtain a split
  exact sequence of graded $\Omega_*^{spin}(pt)$-modules
  \begin{equation}\label{eq:help_describ_K_to_oemge}
    0\to D\to \Omega_*^{spin}(pt)[x]/(Bx-1)\to KO_*(pt)\to 0,
  \end{equation}
  where the middle term is the quotient of the polynomial ring by the ideal
  generated by $(Bx-1)$ where $B$ is the Bott manifold, and where $x$ has
  degree $-8$ and is mapped to the generator of $KO_{-8}(pt)$. Note that
  $\Omega^{spin}_*(pt)[x]/(Bx-1)$ is actually the localization
  $\Omega^{spin}_*(pt)[B^{-1}]$ where we invert $B$. 
The split is determined
  by the inverse of $\tau$ in degrees $0$ through $7$ (where $\tau$ is
  invertible), and by mapping the generator of $KO_8(pt)$ to
  $B$.

  Every element in $\Omega_*^{spin}(pt)[B^{-1}]$ can be written
  (non-uniquely) as $x^k [M]$ for a suitable spin manifold $M$ (because of the
  simple form of the relation, one can multiply every monomial with $B^jx^j$ to
  make any polynomial homogeneous and represent the same element in the
  quotient). Because the image of $x$ in $KO_*(pt)$ is a unit, the kernel $D$
  consists therefore of elements  of the form $[M]x^k$
  with $\hat{A}(M)=0$.

Tensor now \eqref{eq:help_describ_K_to_oemge} 
  with $\Omega_*^{spin}(X)$ to obtain a short split exact sequence of abelian groups
  \begin{multline}\label{eq:describe_ker_of_omega_to_K}
    0\to \Omega^{spin}_*(X)\tensor_{\Omega^{spin}_*(pt)}D\\ \to
    \Omega_*^{spin}(X)[B^{-1}]\to
    \Omega_*^{spin}(X) \tensor_{\Omega^{spin}_*(pt)}   KO_*(pt)\to 0.
\end{multline}
Observe now that, by the universal property of the localization, the $\Omega^{spin}_*(pt)$-module homomorphism
 $\Omega^{spin}_*(X)\to
  \Omega_*(X)\tensor_{\Omega_*^{spin}(pt)}KO_*(pt)\iso KO_*(X)$ 
  factors through
  $\Omega^{spin}_*(X)[B^{-1}]$ because $B$ is mapped to an invertible element
 of $KO_*(pt)$. Every element in the kernel of the projection
  $\Omega^{spin}_*(X)[B^{-1}]\to KO_*(X)$ is a product of $B^{-k}$ with a
  disjoint union of elements $[A_i\times C_i,u_i]$ as in the statement of the
  lemma. By assumption, $[M,v]$ is mapped to such an element in the
  localization $\Omega_*^{spin}(X)[B^{-1}]$. Finally, two elements in such a
  localization are equal if they are bordant (i.e.~equal in
  $\Omega^{spin}_*(X))$ after multiplication with a sufficiently high power of
  $B$. This finishes the proof of the lemma.
\end{proof}

We conclude that, since $l\cdot [M,v] = 0 \in KO_*(B[\Gamma_1\times
\Gamma_2])$, there are product manifolds $A_i\times C_i$ with
 $\hat{A} (C_i)=0$
and with continuous maps $u_i\colon A_i\to
B[\Gamma_1\times\Gamma_2]$ and $n\in\naturals$ such that $l\cdot [M,v]\times B^n$ is
bordant in $\Omega^{spin}_*(B[\Gamma_1\times\Gamma_2])$ to $\sum
[A_i,u_i]\times C_i$. Let $[W,f]$ be the corresponding bordism. Note
that the flat bundles associated to $\lambda_1$ and $\lambda_2$ pull
back from $M$ to each copy of $M\times B^n$ in $\boundary W$ and
extend to all of $W$, restricting on $A_i\times C_i$ to flat bundles
which pull back from $A_i$.

By the multiplicativity of eta invariants (compare \cite{Gilkey}),
\begin{equation*}
\rho_{\lambda_1-\lambda_2}(M\times B^n) =
\rho_{\lambda_1-\lambda_2}(M)\cdot \hat{A}(B^n) =
\rho_{\lambda_1-\lambda_2}(M).
\end{equation*}

On the other hand,
\begin{equation*}
  \rho_{\lambda_1-\lambda_2}(A_i\times C_i) =
  \rho_{\lambda_1-\lambda_2}(A_i)\hat{A}(C_i)=0 .
\end{equation*}

By the  classical  Atiyah-Patodi-Singer index theorem, the
difference of the 
APS-indices of the Dirac operator on $W$ twisted with the flat bundles
associated to $\lambda_1$ and $\lambda_2$ is an integer, which is equal
 to the difference of the
rho-invariants of $l\cdot M\times B^n$ and of $A_i\times C_i$ (indeed, the local terms 
will cancel out). To
conclude
\begin{equation*}
  l\cdot \rho_{\lambda_1-\lambda_2}(M) \in \integers,
\end{equation*}
and this is exactly what we had to prove.

It should be remarked that the calculation of
$KO_*(X)$ in terms of spin-bordism of \cite{Ho-Ho} is very
non-trivial, and is crucially used in our argument.

{\small \bibliographystyle{plain}
\bibliography{rho}
}

\end{document}